\documentclass[11pt]{amsart}
\usepackage[dvipsnames,usenames]{color}
\usepackage{hyperref}
\usepackage{comment}
\usepackage{amsmath}
\usepackage{amsfonts}
\usepackage{amssymb}
\usepackage{amsthm}
\usepackage{amscd}
\usepackage[total={6.6in,8.8in},
top=1in, left=0.9in, includefoot]{geometry}

\newtheorem*{theorem-non}{Theorem}

\newtheorem*{conjecture-non}{Conjecture}

\newtheorem*{lemma-non}{Lemma}
\theoremstyle{remark}

\newtheorem*{remark-non}{Remark}

\def\Z{\mathbb{Z}}

\title{ Cosmetic surgeries and the Poincar\'e homology sphere  }

\author[Tye Lidman]{Tye Lidman}
\thanks{The author was partially supported by NSF grant DMS-1709702 and a Sloan Fellowship.}
\address{Department of Mathematics, North Carolina State University, Raleigh, NC 27607}
\email{tlid@math.ncsu.edu}

\numberwithin{equation}{section}

\begin{document}

\maketitle \vspace{-.1in} 

\begin{abstract}
\vspace{-.25in}  In this short note, we prove that if a knot in the Poincar\'e homology sphere is homotopically essential, then it does not admit any purely cosmetic surgeries. \\
\end{abstract}

Let $K$ be a knot in a three-manifold $Y$.  Recall that a knot $K$ admits  {\em purely cosmetic surgeries} if two different Dehn surgeries on $K$ result in orientation-preserving homeomorphic manifolds.\footnote{We remark that there are many {\em chirally} cosmetic surgeries, i.e. surgeries that result in orientation-reversing homeomorphic manifolds.  (See, for instance, \cite{Mathieu}.)  We will only focus on purely cosmetic surgeries in this note.} 
  Determining when knots have such surgeries is an interesting problem in low-dimensional topology, as it subsumes the knot complement problem and the nugatory crossing conjecture.  Note that the unknot in $Y$ admits purely cosmetic surgeries.  However, the Cosmetic Surgery Conjecture (see \cite[Conjecture 6.1]{Gordon}) asserts that this is the only such knot:

\begin{conjecture-non}
Let $Y$ be a closed, oriented three-manifold and $K$ a knot whose exterior is boundary-irreducible.  If there exist two surgery slopes for $K$ which produce orientation-preserving homeomorphic manifolds, then there is a homeomorphism of the exterior sending one slope to the other.
\end{conjecture-non}  

This problem is an effective testing ground for various Dehn surgery techniques.  Consequently, many restrictions have been established recently using tools from hyperbolic geometry, curve complexes, character varieties, quantum invariants, and Heegaard Floer homology.  For example, building on work of Ozsv\'ath-Szab\'o \cite{HFKQ} and J. Wang \cite{Wang}, Z. Wu showed that if $p/q$ and $p'/q'$ are purely cosmetic non-trivial surgeries on a non-trivial knot in $S^3$, then $p/q$ and $p'/q'$ have opposite sign \cite{Wu}.  This was extended by Ni and Wu \cite{NiWu} to show that $p/q = -p'/q'$, and gave several other constraints on the surgery coefficients and the knot.  We will not use this strengthening, but will use the fact that Wu's original results apply more generally: they hold for non-trivial knots in any integer homology sphere L-space.  One example of such a 3-manifold is the Poincar\'e homology sphere (and conjecturally all other examples are connected sums thereof).  Using this result, we will deduce that ``most'' knots in the Poincar\'e homology sphere cannot admit purely cosmetic surgeries.  

\begin{theorem-non}\label{thm:main}
Let $K$ be a knot in the Poincar\'e homology sphere which is not nullhomotopic.  Then $K$ does not admit any cosmetic surgeries.    
\end{theorem-non}

The rough strategy is as follows.  We use the cosmetic surgeries to build a positive definite cobordism from the Poincar\'e homology sphere to itself.  By choosing $K$ to be homotopically essential, the resulting cobordism has no non-trivial $SU(2)$ representations.  Such a four-manifold cannot exist by Taubes's periodic ends theorem \cite{Taubes}.  We now give the proof.  

\begin{proof}
Without loss of generality, we orient the Poincar\'e homology sphere as the boundary of the positive-definite $E8$ plumbing and denote this three-manifold by $Y$.  First, suppose that $-\infty < p'/q' < 0 < p/q < \infty$ are the two cosmetic surgery slopes for $K$ and denote the result of the surgeries by $Z$.  Then, there exists a positive-definite 2-handle cobordism, $V_{p/q}$, from $Y$ to $Z$ obtained as follows.  Choose a continued fraction expansion for $p/q = [a_1,\ldots,a_n]$ with $a_1 \geq 1$ and $a_2,\ldots, a_n \geq 2$.  Consider the framed link $L_{p/q}$ whose first component is $K$ with framing $a_1$ and remaining components given by a chain of meridians with framings $a_2,\ldots,a_n$.  By slam dunks, we see that framed surgery on this link is $Y_{p/q}(K) = Z$.  The cobordism $V_{p/q}$ is then obtained by attaching 2-handles corresponding to $L_{p/q}$.  The constraints on the $a_i$ guarantee that $V_{p/q}$ is positive-definite.  (For more details, see \cite[Lemma 2.5]{Owens} or \cite{Auckly}.)  Similarly, there is a negative-definite 2-handle cobordism, $V_{p'/q'}$, from $Y$ to $Z$.  By reversing the orientation of $V_{p'/q'}$ and flipping it upside down, we may glue this to $V_{p/q}$ to obtain a positive-definite cobordism, $W$, from $Y$ to itself.  In the case that $p/q$ (respectively $p'/q'$) is $\infty$, we alternatively define $V_{p/q}$ (respectively $V_{p'/q'}$) as $Y \times I$.\footnote{While this case has been included for self-containedness, this special case has already been established using Heegaard Floer homology.  See \cite{Gainullin, Ravelomanana}.}

Next we claim that the fundamental group of $W$ has no non-trivial $SU(2)$ representations.  Note that $W$ is obtained from attaching 2-handles to $Y$, so $\pi_1(W)$ is a quotient of $\pi_1(Y)$.  Since one of these 2-handles is attached along $K$, we have more specifically that $\pi_1(W)$ is a quotient of $\pi_1(Y) / \langle \langle K \rangle \rangle$.  We claim that $\pi_1(Y) / \langle \langle K \rangle \rangle$, and hence $\pi_1(W)$, has only trivial $SU(2)$ representations.  There are two cases.  The first is that $K$ normally generates $\pi_1(Y)$.  In this case, clearly $\pi_1(W)$ is trivial.  The second case is that the normal closure of $K$ is not all of $\pi_1(Y)$.  There is only one proper normal subgroup of the binary icosahedral group, and the quotient by this subgroup is isomorphic to $A_5$.  However, it is well-known that $A_5$ has no non-trivial $SU(2)$ representations.  We conclude in either case that $\pi_1(W)$ has no non-trivial $SU(2)$ representations.
    
Since $W$ is positive-definite and $\pi_1(W)$ has no non-trivial $SU(2)$ representations, we can build an admissible end-periodic manifold $M = E8 \cup W \cup W \cup \ldots$ as in \cite[Definition 1.3]{Taubes}.  Since $E8$ is simply-connected and we have only attached two-handles in order to obtain $M$, this manifold is also simply-connected.  Taubes's periodic ends theorem \cite[Theorem 1.4]{Taubes} implies that the intersection form of $E8$ is diagonalizable over $\Z$, which is a contradiction.  
\end{proof}

\begin{remark-non}
The above arguments can be applied to give some more general results about knots in arbitrary homology spheres.  For example, let $Z$ be a homology sphere that bounds a simply-connected, non-standard, definite four-manifold.  If $K$ normally generates $\pi_1(Z)$, then a positive surgery on $K$ cannot be homeomorphic to negative surgery on any knot in $Z$.  Additionally, $K$ is determined by its complement.  
\end{remark-non}

\bibliographystyle{alpha}
\bibliography{references}

\begin{thebibliography}{Wan06}

\bibitem[Auc97]{Auckly}
David Auckly.
\newblock Surgery numbers of {$3$}-manifolds: a hyperbolic example.
\newblock In {\em Geometric topology ({A}thens, {GA}, 1993)}, volume~2 of {\em
  AMS/IP Stud. Adv. Math.}, pages 21--34. Amer. Math. Soc., Providence, RI,
  1997.

\bibitem[Gai18]{Gainullin}
Fyodor Gainullin.
\newblock Heegaard {F}loer homology and knots determined by their complements.
\newblock {\em Algebr. Geom. Topol.}, 18(1):69--109, 2018.

\bibitem[Gor91]{Gordon}
Cameron~McA. Gordon.
\newblock Dehn surgery on knots.
\newblock In {\em Proceedings of the {I}nternational {C}ongress of
  {M}athematicians, {V}ol. {I}, {II} ({K}yoto, 1990)}, pages 631--642. Math.
  Soc. Japan, Tokyo, 1991.

\bibitem[Mat92]{Mathieu}
Yves Mathieu.
\newblock Closed {$3$}-manifolds unchanged by {D}ehn surgery.
\newblock {\em J. Knot Theory Ramifications}, 1(3):279--296, 1992.

\bibitem[NW15]{NiWu}
Yi~Ni and Zhongtao Wu.
\newblock Cosmetic surgeries on knots in {$S^3$}.
\newblock {\em J. Reine Angew. Math.}, 706:1--17, 2015.

\bibitem[OS11]{HFKQ}
Peter~S. Ozsv\'{a}th and Zolt\'{a}n Szab\'{o}.
\newblock Knot {F}loer homology and rational surgeries.
\newblock {\em Algebr. Geom. Topol.}, 11(1):1--68, 2011.

\bibitem[OS12]{Owens}
Brendan Owens and Sa\v{s}o Strle.
\newblock Dehn surgeries and negative-definite four-manifolds.
\newblock {\em Selecta Math. (N.S.)}, 18(4):839--854, 2012.

\bibitem[Rav15]{Ravelomanana}
Huygens~C. Ravelomanana.
\newblock {Knot Complement Problem for L-space $\mathbb{Z}HS^3$}, 2015.
\newblock arXiv:1505.00239.

\bibitem[Tau87]{Taubes}
Clifford~Henry Taubes.
\newblock Gauge theory on asymptotically periodic {$4$}-manifolds.
\newblock {\em J. Differential Geom.}, 25(3):363--430, 1987.

\bibitem[Wan06]{Wang}
Jiajun Wang.
\newblock Cosmetic surgeries on genus one knots.
\newblock {\em Algebr. Geom. Topol.}, 6:1491--1517, 2006.

\bibitem[Wu11]{Wu}
Zhongtao Wu.
\newblock Cosmetic surgery in {L}-space homology spheres.
\newblock {\em Geom. Topol.}, 15(2):1157--1168, 2011.

\end{thebibliography}

\end{document}